\documentclass{owrart}

\def\d {{\partial}}

\newcommand{\be}{\begin{equation}}
\newcommand{\ee}{\end{equation}}

\def \eps{{\varepsilon}}

\def \d{{\partial}}

\newcommand{{\bv}}{ \bf v}

\newtheorem{Thm}{Theorem}

\begin{document}


\begin{talk}[]{Thierry Paul}
{Rossby waves trapped by quantum mechanics}
{Paul, Thierry}

\noindent
Rossby and Poincar\'e waves appear naturally in the study of large scale oceano-graphy. Poincar\'e waves (PW), of period of the order of a day, are fast dispersive
waves and are due to the  rotation of the Earth through the Coriolis force. Much slower, Rossby waves (RW) are sensitive to the 
variations of the Coriolis parameter,
propagate only eastwards and remain localized for long period of times. We would like here to report on some new results, obtained in collaboration with C. Cheverry, I. Gallagher and L. Saint-Raymond \cite{CGPSR1, CGPSR2, CGPSR3} studying this phenomenon, dispersivity of PW and trapping of RW, as a consequence of the study of  the oceanic waves in a shallow water 
flow subject to strong
wind forcing and rotation,  linearized around a inhomogeneous (non zonal) stationary
profile. The main feature of our results, compared to earlier ones, \cite{G,S, GSR} to quote only very few of them, consists in the fact that we abandon both the betaplane approximation (constant Coriolis force) and the zonal aspect (non dependence w.r.t. the latitude) of the convection term (coupling with the wind).

After some scalings and dimensional homogenizations, the Saint-Venant system of equations for the variations $\eta, u$ near a constant value of the height $\bar h$ and divergence free stationary profile of velocity $\bar u$ takes the form (see \cite{CGPSR2, CGPSR3} for details)
 \begin{equation}
\label{SV2}\scriptsize
\begin{array}{rcl}
  \d_t \eta +  \frac1{\eps}  \nabla \cdot u +  \bar u \cdot \nabla \eta +\eps^2 \nabla \cdot (\eta u)&=&0\,\\
 \d_t u  +      \frac1{\eps^2}  b  u^\perp  +\frac1\eps  \nabla \eta +\bar u \cdot \nabla u + u\cdot \nabla \bar u +\eps^2 u \cdot \nabla u &=&0\,
\end{array}
\end{equation}
where~$b$ is the horizontal component of the Earth rotation vector normalized to one and $\eps^{-1}$ measures the Coriolis force.

The linear version of~(\ref{SV2}) reads (here $D := \frac 1 i \partial$ and $x=(x_1,x_2)\in\mathbb R^2$):
     \begin{equation}\label{linsystcgps}
     \displaystyle \eps^2 i\partial_t{{\bv}} + A(x,\eps D, \eps){\bv} =  0, \,  \qquad 
{\bv} = (v_0,v_1,v_2)=(\eta,u_1,u_2) ,
\end{equation} 
with the linear propagator 

\begin{equation}\label{defA}\scriptsize
A (x,\eps D, \eps):=  i  \left( \begin{array}{ccc}   \eps\bar u\cdot \eps \nabla  &\eps \d_1&\eps \d_2 \\
\eps \d_1 & \eps\bar u\cdot \eps \nabla +\eps^2 \partial_1 \bar u_1  & -b +\eps^2 \partial_2 \bar u_1 \\
\eps \d_2 & b  +\eps^2 \partial_1 \bar u_2&  \eps\bar u\cdot \eps \nabla +\eps^2 \partial_2 \bar u_2 \end{array} \right) \,.
\end{equation}

We will concentrate on \eqref{linsystcgps} with the condition that, essentially, $b$ is increasing at infinity with all derivatives bounded in module by $\vert b\vert$ and only non degenerate critical points. Moreover $\bar u$ will have to be smooth with compact support. 

A simplified version of our main result reads as follows (see \cite{CGPSR3} for details).
\begin{Thm}\label{main}
Under certain microlocalization properties of the initial condition, the solution ${\bv}_\eps(t)={\bv}_\eps(t,.)$ of  \eqref{linsystcgps} decomposes on two Rossby and Poincar\'e vector fields ${\bv}_\eps(t)={\bv}^R_\eps(t)+{\bv}^P_\eps(t)$ satisfying
\begin{itemize}
\item $\forall t>0,\  \forall\Omega$ compact set of $\mathbb R^2$, 
\be
\hskip 1.3cm\Vert{\bv}^P_\eps(t)\Vert_{L^2(\Omega)}=O(\eps^\infty)
\ee
\item  $\exists \Omega$ bounded set of $\mathbb R$ such that, $\forall t>0$ 
\be
\Vert{\bv}^R_\eps(t)\Vert_{L^2(\mathbb R_{x_1}\times(\mathbb R\backslash\Omega)_{x_2})}
=O(\eps^\infty).
\ee
\end{itemize}
\end{Thm}
Theorem \ref{main} shows clearly the different nature of the two type of waves: dispersion for Poincar\'e and confining in $x_2$ for Rossby.
The method of proving Theorem \ref{main} will consist in diagonalizing the ``matrix" $A (x,\eps D, \eps)$. Such diagonalization, if possible, would immediately solve \eqref{linsystcgps} by reducing it to the form 
$\eps^2 i\partial_t{{\bf u}} + D(x,\eps D, \eps){\bf u} =  0$ with $D(x,\eps D, \eps)$ diagonal and solving it component by component. Diagonalizing matrices with operator valued entries is not a simple task, but our next result will show how to achieve it modulo $\eps^\infty$ in the case of matrices with $\eps$-semiclassical type operators entries.

To any (regular enough) function $\mathcal A_\eps\sim\sum_0^\infty\eps^l\mathcal A_l$ on $\mathbb R^{2n}=T^*\mathbb R^n$, possibly matrix valued, we associate the operator $A_\eps$ (densely defined) on $L^2(\mathbb R^n)$ defined by:
\centerline{$f\to A_\eps f,\ (A_\eps f)(x)=\int\mathcal A_\eps  (\frac{x+y}2,\xi)e^{i\frac{\xi(x-y)}\eps}f(y)\frac{dyd\xi}{\eps^n}$.} 

\noindent $\mathcal A_\eps$ is  called the symbol of $A_\eps$ and $A_\eps$ the (Weyl) quantization of $\mathcal A_\eps$.
\newcommand{\mAe}{\mathcal A_\eps}
 \newcommand{\mAo}{\mathcal A_0}

Let $A_\eps$ be such $N\times N$ operator valued matrix of symbol $\mathcal A_\eps\sim\sum_0^\infty\eps^l\mathcal A_l$. We will suppose that 
 $\mAo(x,\xi)$ is Hermitian and therefore is diagonalizable (at each point) by $\mathcal U=\mathcal U(x,\xi)$, $\mathcal U^*\mAo\mathcal U=\mbox{diag}(\lambda_1,\dots,\lambda_N):=\mathcal D$. We will suppose moreover that 
 \be\label{cond}
 \forall (x,\xi),\ \forall i\neq j,\ \vert\lambda_i(x,\xi)-\lambda_i(x,\xi)\vert\geq C>0.
 \ee
 \begin{Thm}[\cite{CGPSR3}]\label{diago}
 There exist $V_\eps$ semiclassical operator and $D_\eps$ diagonal (w.r.t. the $N\times N$ structure) such that
 \[
 V_\eps^{-1}A_\eps V_\eps=D_\eps+O(\eps^\infty)\ \ and \ \ V_\eps^*V_\eps=Id_{L^2(\mathbb R^n, \mathbb C^N)}+O(\eps^\infty)=V_\eps V_\eps^*+O(\eps^\infty).\]
 Moreover $D_\eps=D+\eps D_1+O(\eps^2)$, where $D$ is the Weyl quantization of $\mathcal D$ and 
 $
 D_1$ is the diagonal part of $(\Delta_1-\frac{DI_1+I_1D}2)$ with  ($U$ being the Weyl quantization of $\mathcal U$)
 \be\scriptsize
 \Delta_1=\frac{U^*A_\eps U-D}\eps|_{\eps=0},\ \ \ \ I_1=\frac{U^*U-Id_{L^2(\mathbb R^n, \mathbb C^N)}}\eps|_{\eps=0}.\ee
 \end{Thm}
Let us go back now to the case given by 
 \eqref{defA}. One checks easily that $A (x,\eps D, \eps)$  is of semiclassical type. Its symbol is
\begin{equation}\label{defcalA}
\scriptsize
\mathcal A (x,\xi, \eps)=\left( \begin{array}{ccc}   \eps\bar u\cdot \xi  &\xi_1&\xi_2 \\
\xi_1 & \eps\bar u\cdot \xi +\eps^2 \partial_1 \bar u_1  & -b +\eps^2 \partial_2 \bar u_1 \\
\xi_2 & b  +\eps^2 \partial_1 \bar u_2&  \eps\bar u\cdot \xi +\eps^2 \partial_2 \bar u_2 \end{array} \right)
=\left( \begin{array}{ccc}   0  &\xi_1&\xi_2 \\
\xi_1 & 0  & -b \\
\xi_2 & b  &  0 \end{array} \right)+O(\eps).
\end{equation}
The spectrum of the  leading order $\mathcal A (x,\xi, 0)$ is $\{-\sqrt{\xi^2+b^2(x_2)},0,+\sqrt{\xi^2+b^2(x_2)}\}$. Therefore Condition \eqref{cond} is satisfied only if $\xi^2+b^2(x)\geq C>0$ which correspond to the microlocalization condition in Theorem \ref{main}. Theorem \ref{diago} gives, after a tedious computation, that $A (x,\eps D, \eps)$  is unitary equivalent (modulo $\eps^2$) to the diagonal matrix $\mbox{diag}(T^+,T^R,T^-)$ where $T^\pm$ is the Weyl quantization of $\tau^\pm(x,\xi):=\pm\sqrt{\xi^2+b^2(x_2)}$ and $T^R$ is the quantization of the Rossby Hamiltonian 
$\tau^R(x,\xi):=\eps(\frac{\xi_1b'(x_2)}{\xi^2+b^2(x)}+\bar u(x)\cdot \xi)$.

Under the betaplane approximation, $b(x_2)=\beta x_2$, the Hamiltonians $T^\pm$ are exactly solvable and one shows by hand the dispersive effect for the Poincar\'e waves. In our situation this doesn't work, and because of the $\eps^2$ term in the r.h.s. of \eqref{linsystcgps} the method of characteristics does not apply. A general argument, inherited form quantum mechanics will provide us the solution. First we remark that the Poisson bracket $\{\tau^\pm, x_1\}=\xi_1/\tau^\pm$. This indicates, at a classical level, that $\dot x_1$ has a sign for each Poincar\'e polarization, leading to no return travel. The following theory, due to Eric Mourre, gives the ``quantum" equivalent of this argument.

Let $H$ and $A$ be two self-adjoint operators on a Hilbert space $\mathcal H$ such that: the intersection of the domains of $H$ and $A$ is dense in the domain of $H$, $t \mapsto e^{itA} $ maps~the domain of $H$ to itself and $\sup_{[0,1]}\Vert He^{itA}\varphi\Vert<\infty$ for $\varphi$ in the domain of $H$, and $ i[H,A]$ is bounded from below, closable and the domain of its closure contains the domain of $H$. Finally let us suppose the following

\noindent
\textbf{Positivity condition}: there exist $\theta>0$ and  an open interval~$\Delta$ of~$\mathbb R$ such that if~$E_\Delta$ is the corresponding spectral projection of~$H$, then
\begin{equation} \label{mourreestimate}
E_\Delta i  [H,A] E_\Delta \geq \theta E_\Delta,
\end{equation}
namely $i  [H,A]>0$ on any spectral interval of $H$ contained in $\Delta$.
\begin{Thm}[E. Mourre '80, \cite{JMP}]\label{mourre}
For any integer~$m \in \mathbb N$ and  for any~$\theta ' \in ]0 , \theta [$, there is a constant~$C$ such that
 $$
\|  \chi_- (A-a-\theta't) e^{-iHt} g(H)  \chi_+ (A-a ) \| \leq C  t^{-m}
$$
where~$\chi_\pm$ is the characteristic function of~$\mathbb R^\pm$, $g $ is any smooth compactly supported function in~$\Delta$, and the above bound is uniform in~$a \in \mathbb R$. 
\end{Thm}
In other words, to talk in the quantum langage, if one starts with an initial condition $\varphi$ such that ``$A\geq a$" and the positivity condition \eqref{mourreestimate} holds, after any time $t$ the ``probability" that ``$A\leq \theta't$" is  of order $t^{-m}$. In particular, as $t\to\infty$ the solution $e^{-iHt}\varphi$ escape from any compact spectral region of $A$. 

\vskip 0.3cm
Taking $A=x_1$, Theorem \ref{mourre} gives, after verification that it applies, exactly the ``Poincar\'e" part of Theorem \ref{main}. The ``Rossby part" is given by using the bicharacteristic method and a small computation done in \cite{CGPSR3} which shows that bicharacteristics are trapped in finite regions in the latitude ($x_2$) direction.

Let us mention to finish that the nonlinear terms  can be handled by using a ``$L^\infty$" Gronwall Lemma  and working in some anisotropic and semiclasical Sobolev spaces, so that the solution of \eqref{SV2} is close to the one of \eqref{linsystcgps} as $\eps\to 0$.

\end{talk}

\end{document}